\documentclass[12pt]{amsart}
\usepackage{upref}
\usepackage[latin1]{inputenc}
\usepackage{mathptmx}

\let\oldlabel=\label
\def\prellabel{\marginparsep=1em\marginparwidth=44pt
     \def\label##1{\oldlabel{##1}\ifmmode\else\ifinner\else
          \marginpar{{\footnotesize\ \\ \tt
                     ##1}}\fi\fi}}

\newtheorem{theorem}{Theorem}[section]
\newtheorem{lemma}[theorem]{Lemma}
\newtheorem{corollary}[theorem]{Corollary}
\newtheorem{proposition}[theorem]{Proposition}
\theoremstyle{definition}
\newtheorem{remark}[theorem]{Remark}

\newtheorem{example}[theorem]{Example}

\def\NZQ{\mathbb}               
\def\NN{{\NZQ N }}

\def\ZZ{{\NZQ Z}}

\def\opn#1#2{\def#1{\operatorname{#2}}} 

\def\mm{{\mathfrak m}}
\def\MM{{\mathfrak M}}
\def\pp{{\mathfrak p}}
\def\qq{{\mathfrak q}}
\def\Rees{{\mathcal R}}
\def\DD{{\mathcal D}}
\opn\htt{ht} \opn\cl{cl} \opn\Cl{Cl} \opn\Pic{Pic} \opn\Hom{Hom}
\opn\Ext{Ext} \opn\chara{char} \opn\R{R} \opn\gr{gr}
\opn\Spec{Spec} \opn\rank{rank}

\def\implies{\ifmmode\Rightarrow \else
     \unskip${}\Rightarrow{}$\ignorespaces\fi}

\begin{document}

\title{Canonical modules of Rees algebras}

\author{Winfried Bruns \and Gaetana Restuccia}

\address{Universit\"at Osnabr\"uck,
FB Mathematik/Informatik, 49069 Osnabr\"uck, Germany}

\email{Winfried.Bruns@mathematik.uni-osnabrueck.de}

\address{Università di Messina, Dipartimento di Matematica,
Contrada Papardo, salita Sperone, 31, 98166 Messina}
\email{grest@dipmat.unime.it}

\dedicatory{To Wolmer Vasconcelos on his 65th birthday}


\begin{abstract}
We compute the canonical class of certain Rees algebras. Our
formula generalizes that of Herzog and Vasconcelos.  Its proof
relies on the fact that the formation of the canonical module
commutes with subintersections in important cases. As an
application we treat the classical determinantal ideals and the
corresponding algebras of minors.
\end{abstract}

\maketitle

A considerable part of Wolmer Vasconcelos' work has been devoted
to Rees algebras, in particular to their divisorial structure and
the computation of the canonical module (see \cite{HSV1, HSV2, HV,
MV, Va}). In this paper we give a generalization of the formula of
Herzog and Vasconcelos \cite{HV} who have computed the canonical
module under more special assumptions. We show that
$$
[\omega_\Rees]=[I\Rees]+\sum_{i=1}^t (1-\htt\pp_i)[P_i]
$$
for ideals $I$ in regular domains $R$ essentially finite over a
field for which the Rees algebra $\Rees=\Rees(I)$ is normal
Cohen--Macaulay and whose powers have a special primary
decomposition. In the formula $P_1,\dots,P_t$ are the divisorial
prime ideals containing $I\Rees$, and $\pp_i=P_i\cap R$. The
condition on the primary decomposition can be expressed
equivalently by the requirement that the restriction of the Rees
valuation $v_{P_i}$ to $R$ coincides with the $\pp_i$-adic
valuation.

It is not difficult to derive a criterion for the Gorenstein
property of the Rees algebra, the extended Rees algebra and the
associated graded ring from the formula above.

While our hypotheses are far from the most general case, which may
very well be intractable, the formula covers many interesting
ideals, for example the classical determinantal ideals. As an
application we can therefore compute the canonical classes of
their Rees algebras, and also those of the algebras generated by
minors. The algebras generated by minors can be identified with
fiber cones of the determinantal ideals, and therefore are
accessible via the Rees algebra.

This paper has been inspired by the work of Bruns and Conca
\cite{BC2} where the case of the determinantal ideals of generic
matrices and Hankel matrices has been treated via initial ideals.

The formula above can be proved with localization arguments only
if there is no containment relation between the prime ideals
$\pp_i$. It was therefore necessary to investigate the behavior of
the canonical module under subintersections, and, as we will show
for normal algebras essentially of finite type over a field, its
formation does indeed commute with taking a subintersection.

\section{Divisor class group and valuations on a Rees algebra}

Let $R$ be a normal Noetherian domain, and $I$ an ideal in $R$ for
which the Rees algebra
$$
\Rees=\Rees(I)=\Rees(R,I)=\bigoplus_{k=0}^\infty I^kT^k\subset
R[T]
$$
is a normal domain. The normality of the Rees algebra is
equivalent to the integral closedness of all powers $I^k$ (see
Ribenboim \cite{Ri}). We will assume in the following that $\htt
I\ge 2$.

The divisor class group of $\Rees$ has been investigated in
several articles of which Vasconcelos is a co-author (\cite{HSV1,
HSV2, HV, MV}). The main result, first proved by Simis and Trung
\cite{ST}, describes $\Cl(\Rees)$ as follows:

\begin{theorem}\label{ClRees}
There is an exact sequence
\begin{equation}
0\to\ZZ^t\to\Cl(\Rees)\to \Cl(R)\to 0\label{ClSeq}
\end{equation}
where a basis of $\ZZ^t$ is given by the classes $[P_i]$ of the
minimal prime ideals $P_1,\dots,P_t$ of the divisorial ideal
$I\Rees$.
\end{theorem}

In the following we will denote the subgroup
$\ZZ[P_1]+\dots+\ZZ[P_t]$ simply by $\ZZ^t$.

The most informative proof of this theorem uses the following
lemma, in which $\Spec^1$ denotes the set of divisorial prime
ideals.

\begin{lemma}\label{ReesVal}
Let $I$ be an ideal of height at least
$2$, and, with the notation introduced, $V_i=\Rees_{P_i}$; then
\begin{equation}
\Rees=R[T]\cap V_1\cap\dots\cap V_t,\label{ReesInt}
\end{equation}
and, moreover,
\begin{equation}
R[T]=\bigcap \{\Rees_P: P\in \Spec^1(\Rees),\ P\neq
P_1,\dots,P_t\}.\label{SubInt}
\end{equation}
\end{lemma}

\begin{proof}
Since $\Rees$ is a Krull domain, it is the intersection of its
localizations $\Rees_P$, where $P$ runs through $\Spec^1(\Rees)$.
For \eqref{ReesInt} we have to show that either $\Rees_P=R[T]_Q$
for a divisorial prime ideal $Q$ of $R[T]$, or $P=P_i$ for some
$i$.

If $P\supset I\Rees$, then $P=P_i$ for some $i$. Otherwise there
exists $x\in I$, $x\notin P$. Then $\Rees_P$ is a localization of
$\Rees[x^{-1}]$, and, since $\Rees[x^{-1}]=R[T,x^{-1}]$, it is
also a localization of $R[T]$ with respect to some divisorial
prime ideal.

Equation \eqref{SubInt} follows from \cite[Lemma 2.1]{BC2} because
the extensions $P_iR[T]$, containing $IR[T]$, have height at least
$2$. (One should beware of considering \eqref{SubInt} as a trivial
consequence of \eqref{ReesInt}; for $R=K[X]$, $K$ a field, and
$I=RX$, the intersection in \eqref{SubInt} is $R[X^{-1},T]$.)
\end{proof}

Equation \eqref{SubInt} says that $R[T]$ is a subintersection of
$\Rees$ in the terminology of Fossum \cite{Fo}, and so Nagata's
theorem \cite[Thm.\ 7.1]{Fo} yields the exact sequence
\eqref{ClSeq} with $\Cl(R[T])$ in place of $\Cl(R)$, once one has
shown that the classes $[P_1],\dots,[P_t]$ are linearly
independent. A proof of the linear independence can be found in
Morey and Vasconcelos \cite{MV}. Finally one uses that
$\Cl(R[T])=\Cl(R)$.

We want to complement Theorem \ref{ClRees} by a description of the
Picard group of the Rees algebra. It is the group of isomorphism
classes of projective rank $1$ modules, and therefore a natural
subgroup of the class group. One has $\Pic(R)=\Cl(R)$ if and only
if $R$ is locally factorial (i.~e.\ all localizations of $R$ at
maximal ideals are factorial).

\begin{proposition}\label{PicRees}
The natural map $\Pic(R)\to\Pic(\Rees)$ is an isomorphism.

If $R$ is locally factorial, then the natural map
$\Cl(R)=\Pic(R)\to\Pic(\Rees)\subset\Cl(\Rees)$ splits the exact
sequence \eqref{ClSeq}. In particular
$\Cl(\Rees)=\ZZ^t\oplus\Pic(\Rees)$.
\end{proposition}

\begin{proof}
Let us first show that $\Pic(\Rees)=0$ if $R$ is a local ring with
maximal ideal $\mm$. In this case $\Rees$ is a graded ring whose
homogeneous non-units generate the maximal ideal
$\MM=\mm\oplus\bigoplus_{k=1}^\infty I^kT^k$. Moreover, each
divisorial ideal of $\Rees$, especially every rank $1$ projective
module, is isomorphic to a graded ideal. A graded projective
module is free, since its localization with respect to $\MM$ is
free (see \cite[1.5.15]{BH}).

We now turn to the general case. The natural map
$\Pic(R)\to\Pic(R[T])$, induced by ring extension, is an
isomorphism, since $R$ is normal. It factors through
$\Pic(\Rees)$. Therefore the map $\Pic(\Rees)\to\Pic(R[T])$ is
surjective. In order to show that it is injective, we have to
verify that $\Pic(\Rees)\cap\ZZ^t=0$. This follows from the fact
that each nonzero divisor class $[C]$ in $\ZZ^t$ survives in at
least one Rees algebra $R_\pp[I_\pp T]$. By the first part of the
proof this is impossible for the class of a projective rank one
module. In fact, if the coefficient of $[C]$ with respect to
$[P_i]$ is nonzero, we choose $\pp=\pp_i$.

The second statement is now obvious.
\end{proof}

The valuations $v_{P_i}$ on the quotient field $Q(\Rees)=Q(R[T])$
are called the \emph{Rees valuations} of $I$ (cf.\ MacAdam
\cite[Ch.\ XI]{MA}). If $v$ is a valuation on a domain $R$ such
that $v(x)\ge 0$ for all $x\in R$, then the \emph{center} of $v$
is the prime ideal $\{x:v(x)>0\}$.

The representation \eqref{ReesInt} of $\Rees$ can be translated
into a description of the powers of $I$ as an intersection of
valuation ideals.

\begin{proposition}\label{IdInt}
Let $v_i$ be the Rees valuation on the quotient field $Q(R[T])$
associated with $\Rees_{P_i}$, $i=1,\dots,t$, and set
$J_i(j)=\{x\in R:v_i(x)\ge j\}$.  Then
\begin{equation}
I^k=\bigcap_{i=1}^t J_i(kd_i),\qquad d_i=-v_i(T).\label{IdIntE}
\end{equation}
The intersection is irredundant for $k\gg 0$. Moreover,
$$
I\Rees=\bigcap_{i=1}^t P_i^{(d_i)}.
$$
\end{proposition}

\begin{proof}
We consider equation \eqref{ReesInt} in each $T$-degree. Then it
says
$$
I^kT^k=\{aT^k:a\in R,\ v_i(a)\ge -v_i(T^k),\ i=1,\dots,t\},
$$
and this is evidently equivalent to equation \eqref{IdIntE}.

To see that the intersection is irredundant for $k\gg 0$, we use
that the representation \eqref{ReesInt} is irredundant: there
exists $y_i\in R$ and $n_i\in\NN$ such that $v_i(y_iT^{n_i})<0$,
but $v_j(y_iT^{n_i})\ge 0$ for $j\neq i$. Moreover,
$IT\Rees=\Rees\cap TR[T]$ is a divisorial prime ideal different
from $P_i$ for all $i$. Thus we can find $x_i\in I$ such that
$v_i(x_iT)=0$. It follows that $v_i(x_i^my_iT^{n_i+m})<0$ for all
$m\ge 0$, but $v_j(x_i^my_iT^{n_i+m})\ge0$ for $j\neq i$. To sum
up: the representation is irredundant for $n\ge \max_i n_i$.

For the second formula we note that $IT\Rees$ is a prime ideal
different from $P_1\dots,P_t$, the divisorial prime ideals
containing $I\Rees$. So $0=v_i(IT\Rees)=v_i(T)+v_i(I\Rees)$.
Therefore $v_i(I\Rees)=-v_i(T)$.
\end{proof}

Let $\pp_i=P_i\cap R=J_i(1)$. Then one sees immediately that
$J_i(j)$ is $\pp_i$-primary for all $j$, and thus \eqref{IdIntE}
yields a primary decomposition of $I^k$ for all $k$. But even if
the intersection is irredundant, it need not be an irredundant
primary decomposition in the usual sense, since it may very well
happen that $\pp_i=\pp_j$ for $i\neq j$.

In the next section we want to compute the class of the canonical
module of $\Rees$ in certain cases, in which we can identify the
ideals $J_i(j)$.

Suppose that $R$ is a regular domain. Then each prime ideal $\pp$
of $R$ defines a discrete valuation on $Q(R)$ as follows. First we
replace $R$ by $R_\pp$, and may assume that $R$ is local with
maximal ideal $\pp$. Now  we set
$v_\pp(x)=\max\{i:x\in\pp^{(i)}\}$ for each $x\in R$, $x\neq 0$,
and extend this function naturally to $Q(R)$ (with
$v_\pp(0)=\infty$). That $v_\pp$ is indeed a valuation follows
from the fact that the associated graded ring of the filtration
$(\pp^i)$ is a polynomial ring over the field $R/\pp$ and
therefore an integral domain. It guarantees that
$v_\pp(xy)=v_\pp(x)+v_\pp(y)$. A similar argument shows that the
symbolic powers of prime ideals in regular domains are integrally
closed.

One says that a valuation on the polynomial ring $R[T]$ (where $R$
may be an arbitrary domain) is \emph{graded} if
$$
v(f)=\min_i v(a_iT^i)
$$
for all polynomials $f=\sum a_iT^i$, $a_i\in R$. Every valuation
on $R$ can be extended to a graded valuation on $R[T]$. One can
freely choose $v(T)$ and then use the previous equation to define
the extension of $v$ (see Bourbaki \cite[Ch.\ VI, \S10, no.\ 1,
lemme 1]{Bo}).

The Rees valuations of $I$ on $R[T]$ are graded. In fact, if $S$
is a normal graded subalgebra of $R[T]$ with $Q(S)=Q(R[T])$, then
the valuations associated with graded divisorial ideals of $S$ are
graded, since all symbolic powers of graded prime ideals are
graded as well. Moreover, the associated prime ideals of the
graded ideal $I\Rees$ are graded.

\begin{proposition}\label{PriDec}
Let $R$ be a regular ring and $I$ an ideal of height $\ge 2$.
\begin{itemize}
\item[(a)] Then the following are equivalent:
\begin{itemize}
\item[(i)] $\Rees(I)$ is normal, and for each minimal prime ideal
$P$ of $I\Rees$ the Rees valuation $v_P$ restricts on $R$ to the
valuation $v_\pp$, $\pp=P\cap R$;

\item[(ii)] there exist prime ideals $\pp_i,\dots,\pp_u$ in $R$
and $d_1,\dots,d_u\in\NN$ such that $I^k=\bigcap_{i=1}^u
\pp_i^{(d_ik)}$ for all $k$.
\end{itemize}
\item[(b)] Moreover, if \emph{(i)} holds, and $P_1,\dots,P_t$ are
the minimal prime ideals of $I\Rees$, then one can choose
$\pp_i=P_i\cap R$, $d_i=-v_i(T)$, and the intersection in
\emph{(ii)} is irredundant for $k\gg 0$.

\item[(c)] Conversely, if there exists $k_i$ for each
$i=1,\dots,u$ such that $\pp_i^{(d_ik)}$ cannot be omitted in the
representation of $I^{k_i}$ in \emph{(ii)}, then the graded
extensions of the $v_i$ to $R[T]$ with $v_i(T)=-d_i$ are the Rees
valuations of $I$ on $Q(R[T])$.
\end{itemize}
\end{proposition}

\begin{proof}
(a) The implication (i)\implies(ii) has been proved above,
together with the description of the $\pp_i$ in terms of the prime
ideals $P_i$.

Suppose now that (ii) holds. Then all powers of $I$ are integrally
closed, since the symbolic powers of prime ideals in regular rings
are integrally closed. Let $v_i$ be the valuation on $R[T]$ that
we obtain as the graded extension of $v_{\pp_i}$ to $R[T]$ with
$v_i(T)=-d_i$. The representation of $I^k$, $k\in\NN$, can
immediately be translated into the description of $\Rees(I)$ as
the intersection of $R[T]$ with the discrete valuation rings $V_i$
associated with the valuations $v_i$:
\begin{equation}\label{IntRees1}
\Rees=R[T]\cap V_1\cap\dots\cap V_u.
\end{equation}
Let $P_i$ be the center of $v_i$ in $\Rees$. Then $I\Rees=\bigcap
P_i^{(d_i)}$. If this representation is not irredundant, we can
shorten it to an irredundant primary decomposition, which is
unique since $I\Rees$ is a divisorial ideal. As seen above, we can
shorten the representation \eqref{IntRees1} accordingly. Thus we
may assume that $P_1,\dots,P_u$ are the minimal primes of
$I\Rees$. That the associated valuations satisfy the condition in
(i), follows from their construction.

(b) is only a restatement of Proposition \ref{IdInt} under the
special hypothesis made in (i).

(c) The condition guarantees that none of the $V_i$ can be omitted
in \eqref{IntRees1}, and the rest has been proved above.
\end{proof}

\begin{remark}\label{General}
With the appropriate modifications, the results in this section
remain true if one considers an arbitrary ideal of height $\ge 2$
in a normal (or regular) ring and replaces the (in general
non-normal) Rees algebra by its integral closure,
$$
\overline\Rees=\bigoplus_{k=0}^\infty \overline{I^k}T^k
$$
where $\overline{I^k}$ is the integral closure of $I^k$. (The
ideal $I\Rees$ must be replaced by $\bigoplus
\overline{I^{k+1}}T^k$.)
\end{remark}

\section{The canonical class}

As in the previous section, we assume that $I$ is an ideal in the
normal domain $R$ with a normal Rees algebra $\Rees=\Rees(I)$.
Suppose further that $R$ is a Cohen--Macaulay residue class ring
of a Gorenstein ring, and $\Rees$ is also Cohen--Macaulay. Then
$\Rees$ has a canonical module $\omega_\Rees$ (see Bruns and
Herzog \cite{BH}). The canonical module is (isomorphic to) a
divisorial ideal, and this allows us to find $\omega_\Rees$ by
divisorial computations.

The canonical module is unique only up to tensor product with a
projective rank $1$ module. In other words, only its residue class
modulo the Picard group $\Pic(R)\subset\Cl(R)$ is unique.

Let us assume for the moment (and without essential restriction
for the theorem to be proved) that $R$ is factorial. Because of
the exact sequence \eqref{ClSeq} and since $R$ (and, along with
it, $R[T]$) is factorial, the class of $\omega_\Rees$ is a linear
combination of the classes $[P_i]$,
$$
[\omega_\Rees]=w_1[P_1]+\dots+w_t[P_t].
$$
We have to determine the coefficients $w_i$.

Since the behavior of the class group and the canonical module
under localization is easily controlled, we can first replace $R$
by $R_{\pp_i}$ and $I$ by $IR_{\pp_i}$ in order to compute $w_i$
(as above, $\pp_i=P_i\cap R$). However, in general this
localization does not strip off all the other components
$\ZZ[P_j]$: those with $\pp_j\subset\pp_i$ survive.

Therefore we need a finer instrument to isolate $w_i$: we pass to
the subintersection $R[T]\cap V_i$ (after the localization). Then
we must
\begin{itemize}
\item[(i)] determine the structure of $R[T]\cap V_i$ and find its
canonical class, and

\item[(ii)] show that the canonical module is preserved under
subintersection.
\end{itemize}

Let us first turn to Problem (ii). We cannot present a solution in
complete generality. However, the next theorem should cover many
interesting applications.

\begin{theorem}\label{CanSubInt}
Let $K$ be a field, $R$ a normal Cohen--Macaulay $K$-algebra
essentially of finite type over $K$, and $Y\subset \Spec^1(R)$.
Suppose that the subintersection $S=\bigcap_{\pp\in Y}R_\pp$ is
again essentially of finite type over $K$ and Cohen--Macaulay.

Then the canonical module of $S$ is $(\omega_R\otimes_R
S)^{\dagger\dagger}$, where $\vphantom{(}^\dagger$ denotes the
functor $\Hom_S(\_\,,S)$. In other words, the canonical class of
$S$ is the image of $\omega_R$ under the natural map
$\Cl(R)\to\Cl(S)$.
\end{theorem}

\begin{proof}
For technical simplicity let us first assume that $K$ is a perfect
field, and let $\Omega_{R/K}$ be the module of K\"ahler
differentials of $R$. It has been proved by Kunz \cite{Ku1} that
the canonical module is given by the regular differential
$r$-forms $\R_K^r(R)$, $r=\dim R$, and Platte and Storch \cite{PS}
have noticed that $\R_K^r(R)$ can be identified with the
$R$-bidual $\bigl(\bigwedge^r \Omega_{R/K}\bigr)^{**}$. The same
applies to $S$. (A proof will be given below.)

The extension $\phi:R\to S$ gives rise to an $R$-linear map
$d\phi:\Omega_{R/K}\to\Omega_{S/K}$, and $d\phi$ induces a natural
$S$-linear map
$$
\psi: \biggl(\bigwedge^r \Omega_{R/K}\biggr)\otimes_R S \to
\bigwedge^r \Omega_{S/K}.
$$
Let $\qq$ be a height $1$ prime ideal of $S$. Then
$S_\qq=R_{\qq\cap R}$, and therefore $\psi\otimes_S S_\qq$ is an
isomorphism. It follows that the $S$-bidual extension
$\psi^{\dagger\dagger}$ is an isomorphism at all height $1$ prime
ideals $\qq$ of $S$. Since the $S$-biduals are reflexive,
$\psi^{\dagger\dagger}$ is an isomorphism itself.

The second statement about the divisor classes follows
immediately, since $(I\otimes S)^{\dagger\dagger}$ is exactly the
divisorial ideal of $S$ to which a divisorial ideal $I$ of $R$
extends; see \cite{Fo}.

If $K$ is not perfect (this can happen only in characteristic
$p>0$), one replaces $K$ by a subfield $K_0$ with $[K:K_0]<\infty$
that is admissible in the sense of \cite[6.23]{Ku2} for $R$, $S$
and regular $K$-algebras $A$ and $B$ essentially of finite type
for which there exist presentations $R=A/I$ and $S=B/J$. All the
algebras involved are then essentially of finite type over $K_0$,
too.

Let us now show that $\bigl(\bigwedge^r \Omega_{R/K_0}\bigr)^{**}$
is the canonical module of $R$. To this end we let $[M]$ denote
the divisor class of a finitely generated $R$-module: $[M]$ is the
isomorphism class of $\bigl(\bigwedge^n M\bigr)^{**}$ where
$n=\rank M$.

We use the complex
$$
0\to I/I^2\to \Omega_{A/K_0}\otimes_A R\to \Omega_{R/K_0}\to 0
$$
that is exact at the right and in the middle and exact in $I/I^2$
at all prime ideals $\pp$ for which $R_\pp$ is regular. Moreover
$I/I^2$ is free at such $\pp$ of rank $c=\htt I$ (see \cite{Ku2}).
By divisorial calculation,
$$
[I/I^2]=- [\Omega_{R/K_0}],
$$
since $\Omega_{A/K_0}\otimes_A R$ is a free module, and so has
class $0$. Finally, by Herzog and Vasconcelos \cite[Lemma]{HV},
$[\omega_R]=-[I/I^2]$.
\end{proof}

Before we state and prove the main result, let us single out a
very special case.

\begin{proposition}\label{Regu}
Let $R$ be a regular local ring with maximal ideal $\mm$. Then the
Rees algebra $\Rees_k=\Rees(\mm^k)$ is normal and Cohen--Macaulay.
Its canonical module is unique up to isomorphism and has class
$(k-\dim R+1)[P_k]=[\mm^k\Rees_k]-(\dim R-1)[P_k]$ where
$P_k=\mm\Rees_k$ is the only divisorial ideal of $\Rees_k$
containing $\mm^k\Rees_k$.
\end{proposition}

\begin{proof}
By Proposition \ref{PicRees} one has $\Pic(\Rees_k)=0$, hence the
uniqueness of the canonical module.

For the equation $P_k=\mm\Rees_k$ it is enough to show that
$\mm\Rees_k$ is a prime ideal of height $1$. This follows
immediately from the fact that $\Rees_k/\mm\Rees_k$ is the $k$th
Veronese subring of the associated graded ring of $R$ with respect
to $\mm$, a polynomial ring over the field $R/\mm$.

Let $x_1,\dots,x_r$, $r=\dim R$, be a regular system of parameters
and set
$$
J_k=\bigl(x_1\cdots x_rTR[T]\bigr)\cap P_k.
$$
We know that $[P_k]$ generates $\Cl(\Rees_k)$, and $TR[T]\cap
\Rees_k=\mm^k T\Rees_k\cong \mm^k\Rees_k$. Furthermore
$(x_iR[T])\cap\Rees_k$ has divisor class $-[P_k]$. In fact
$x_i\Rees_k=x_iR[T] \cap P_k$, since $x_iR[T]\cap\Rees_k$ and
$P_k$ are the only divisorial prime ideals of $\Rees_k$ containing
the prime element $x_i$ of $R[T]$ (see Lemma \ref{ReesVal});
moreover, $x_i$ obviously has value $1$ under the corresponding
valuations.

This shows that $J_k$ has the class given in the theorem, and it
only remains to prove that it is the canonical module. In the case
$k=1$ this follows immediately from Herzog and Vasconcelos
\cite{HV}. Before going to general $k$, we want to convince
ourselves that $J_1$ is the graded canonical module of $\Rees_1$
in the sense of \cite[3.6]{BH}. In fact, there is an exact
sequence
$$
0\to\Hom(\Rees_1,J_1)\to\Hom(\mm
T\Rees_1,J_1)\to\Ext_1^\Rees(R,J_1)\to 0.
$$
The graded canonical module of $\Rees_1$ is of the form $J_1(-s)$
for some integer $s$. Exactly for the right choice of $s$,
$\Ext_1^\Rees(R,J_1(-s))$ is the graded canonical module
$\omega_R=R$ of $R$ (with the trivial grading $R_0=R$). We have
only to check that this holds with $s=0$. But $\Hom(\Rees_1,J_1)$
has only components of positive degree, and the degree $0$
component of $\Hom(\mm T\Rees_1,J_1)$ is $x_1\cdots x_rR\cong R$,
as desired.

For general $k$ we have
$$
\omega_{\Rees_k}=\omega_{\Rees(\mm)}\cap\Rees_k=J_1\cap
\Rees_k=J_k
$$
by Goto and Watanabe \cite[3.6.21]{BH} since $\Rees_k$ is the
$k$th Veronese subalgebra of $\Rees_1$ with respect to the
$T$-grading.
\end{proof}

We can now prove the main result.

\begin{theorem}\label{main}
Let $K$ be a field, $R$ be a regular domain essentially of finite
type over $K$, and $I$ an ideal with a Cohen--Macaulay normal Rees
algebra $\Rees=\Rees(I)$. Let $P_1\dots,P_t$ be the divisorial
prime ideals of $I\Rees$, $I\Rees=\bigcap_{i=1}^t P_i^{(d_i)}$,
and suppose that $v_{P_i}|R=v_{\pp_i}$ for $i=1,\dots,t$, with
$\pp_i=P_i\cap R$. Then a module of class
$$
\sum_{i=1}^t (d_i+1-\htt\pp_i)[P_i]= [I\Rees]+\sum_{i=1}^t
(1-\htt\pp_i)[P_i]
$$
is a canonical module $\omega_\Rees$ of $\Rees$. Moreover, $\Rees$
is Gorenstein if and only if $d_i=\htt \pp_i-1$ for all
$i=1,\dots,t$.
\end{theorem}

\begin{proof}
Let $C$ be a module of the class given in the theorem. It is
enough to show that each of its localizations $C_\MM$ with respect
to maximal ideals $\MM$ of $\Rees$ is a canonical module of
$\Rees_\MM$. Such a localization $\Rees_\MM$ is a localization of
$\R_\mm$ with $\mm=R\cap\MM$. Since the definition of $[C]$
commutes with localization in $R$ (in fact, primary decomposition
commutes with such localizations), we may assume that $R$ is
factorial. Then $\Cl(\Rees)=\ZZ^t$, $\Pic(\Rees)=0$ (by
Proposition \ref{PicRees}), and we have a unique isomorphism class
$$
[\omega_\Rees]=w_1[P_1]+\dots+w_t[P_t].
$$
for the canonical module of $\Rees$. It is enough to determine,
say, $w_1$. We localize $R$ with respect to $\pp_1$, and may then
assume that $R$ is regular local with maximal ideal $\mm=\pp_1$.
In the next step we pass to the subintersection $S=R[T]\cap V_1$.
But this subintersection is exactly $\Rees(\mm^{d_1})$, as follows
from Propositions \ref{IdInt} and \ref{PriDec}. (Since $\mm$ is a
maximal ideal, its ordinary and symbolic powers coincide.)

According to Theorem \ref{CanSubInt} the formation of the
canonical class commutes with subintersection, and so $w_1$ is the
coefficient of the canonical module of $\Rees(\mm^{d_1})$ with
respect to the extension of $P_1$. By Proposition \ref{Regu} this
coefficient is $d_1+1-\dim R$, as desired.

Because of the splitting $\Cl(\Rees)=\ZZ^t\oplus \Pic(\Rees)$, the
class given in the theorem is the $\ZZ^t$-component of any
canonical module of $\Rees$. Therefore $\Rees$ is Gorenstein if
and only if the $\ZZ^t$-component vanishes.
\end{proof}

Often it is useful to know the graded canonical module of $\Rees$,
as we have seen in the proof of Proposition \ref{Regu}.

\begin{corollary}\label{Graded}
With the hypotheses of Theorem \ref{main} suppose we can find an
element $x\in R$ such that $v_{P_i}(x)=\htt \pp_i$ for all $i$.
\begin{itemize}
\item[(a)] Then
$$
\omega_\Rees=xTR[T]\cap P_1\cap\dots\cap P_t
$$
is the graded canonical module of $\Rees$ (with respect to the
grading by $T$).

\item[(b)] Suppose that $K$ is infinite, $R$ is the polynomial ring
over $K$ in $n$ variables, graded by total degree, $I$ is graded
and $x$ is homogeneous, and choose $x'\in R$ such that $\deg
x'=n-\deg x$ and $x'\notin \pp_i$ for $i=1,\dots,t$. Then
$xx'TR[T]\cap P_1\cap\dots\cap P_t$ is the bigraded canonical
module with respect to the natural bigrading on $\Rees$.
\end{itemize}
\end{corollary}

This follows as in the special case considered in Proposition
\ref{Regu}. Note that the element $x'$ needed for (b) can always
be found by prime avoidance. (If the ideal $(X_1,\dots,X_n)$ is
among the $\pp_i$, we choose $x'=1$.)

An analogous statement as in (b) holds for monomial ideals. Then
we can choose $x=X_1\cdots X_n$ and obtain the multigraded
canonical module of $\Rees$. This can be proved directly from the
theorem of Danilov and Stanley describing the canonical module of
a normal semigroup ring (see \cite[6.3.5]{BH}).

\begin{example}\label{235}
The following example shows that the condition on the Rees
valuations in Theorem \ref{main} is crucial. Let $I$ be the
integral closure of the ideal $(X^2,Y^3,Z^5)\subset K[X,Y,Z]$, $K$
a field. It has been noticed by Reid, Roberts, and Vitulli
\cite{RRV} (and can easily be checked by {\tt normaliz} \cite{BK})
that $\Rees=\Rees(I)$ is normal. A $K$-basis of $\Rees$ is given
by all monomials $X^aY^bZ^cT^d$ where $15a+10b+6c-30d\ge 0$. In
other words, $\Rees=R[T]\cap V_1$, where the valuation defining
$V_1$ is the multigraded extension of the function that takes the
values $v_1(X)=15$, $v_1(Y)=10$, $v_1(Z)=6$ and $v_1(T)=-30$. It
follows that $[I\Rees]=30[P_1]$. By the theorem of Danilov and
Stanley one has $\omega_\Rees=XYZTR[T]\cap P_1$. Arguing as in the
proof of Proposition \ref{Regu} one obtains
$[\omega_\Rees]=(-15-10-6+30)[P_1]+[P_1]=0$. So $\Rees$ is a
Gorenstein ring. Its canonical module is the principal ideal
generated by $XYZT$.
\end{example}

\begin{remark}\label{Monom}
(a) The hypotheses of the theorem can be weakened. If we
\emph{define} the canonical module via K\"ahler differentials (the
description used in the proof of Theorem \ref{CanSubInt}), then
the hypothesis that the Rees algebra is Cohen--Macaulay is no
longer necessary. The proof of the theorem shows that the
canonical module has class $[\Omega_K(R)\otimes\Rees]+\sum_{i=1}^t
(d_i+1-\htt\pp_i)[P_i]$.

(b) Instead of requiring that $\Rees(I)$ is normal, one could
consider the normalization of $\Rees(I)$. As indicated in Remark
\ref{General}, $I\Rees$ has then to be replaced by $\bigoplus
\overline{I^{k+1}}T^k$.

(c) One can generalize Theorem \ref{main} in such a way that
Example \ref{235} is covered. The first part of its proof, namely
the isolation of each $w_i$, does not use the hypothesis on the
Rees valuations. Therefore, as soon as one can compute the
canonical module of $R[T]\cap V_i$ for each $i$, a generalization
is possible. A suitable hypothesis generalizing the condition
$v_{P_i}|R=v_{\pp_i}$ is the following: there exists a regular
system of parameters $x_1,\dots,x_m$ of $R_{\pp_i}$ such that each
of the ideals $\{x\in R_{\pp_i}:v_{P_i}(x)\ge k\}$ is generated by
monomials in $x_1,\dots,x_m$. Then one can replace $\htt\pp_i$ in
the theorem by $v_{P_i}(x_1\cdots x_m)=\sum_{j=1}^m v_{P_i}(x_j)$.

However, there exist valuations that do not allow such a
``monomialization''. A counterexample was communicated by D.
Cutkosky.
\end{remark}

In view of Theorem \ref{main} it is not difficult to decide when
the extended Rees algebra or the associated graded ring are
Gorenstein. (We are grateful to S. Goto for suggesting the
inclusion of the corollary.)

\begin{corollary}\label{extRees}
Suppose that $R$ and $I$ satisfy the hypothesis of Theorem
\ref{main}. Then the extended Rees algebra
$\widehat\Rees=\Rees[T^{-1}]$ or, equivalently, the associated
ring $\gr_I(R)=\Rees/I\Rees$, is Gorenstein if and only if there
exist $c_i\in \NN$, $i=1,\dots,t$, such that
\begin{itemize}
\item[(i)] $c_id_i=\htt \pp_i-1$ for all $i=1,\dots,t$, and
\item[(ii)] $c_i=c_j$ whenever there exists a maximal ideal $\mm$
of $R$ with $\pp_i,\pp_j\subset\mm$.
\end{itemize}
\end{corollary}

\begin{proof}
The Cohen--Macaulay property of $\Rees$ is inherited by
$\widehat\Rees$ and $\gr_i(R)$ (if $R$ is Cohen-Macaulay), as is
well-known.

The Gorenstein property of $\widehat\Rees$ is local with respect
to $\Spec R$, and the same holds for the associated graded ring.
We can therefore assume that $R$ is local with maximal ideal
$\mm$. Then $\Pic(\widehat\Rees)=0$, as follows by the same
argument as in the proof of Proposition \ref{PicRees}. Moreover,
$\widehat\Rees$ is Gorenstein if and only $\gr_I(R)$ is so, since
the latter is the residue class ring by the homogeneous regular
element $T^{-1}$.

The extended Rees algebra $\widehat\Rees$ is a subintersection of
$\Rees$, namely
$$
\widehat\Rees=\bigcap\{\Rees_Q: Q\in\Spec^1(\Rees), Q\neq
IT\Rees\},
$$
and we can apply Theorem \ref{CanSubInt}. (Note that
$IT\Rees=TR[T]\cap\Rees$.) Thus its divisor class group is
$\Cl(\Rees)/\ZZ[IT\Rees]$ (this was noticed in \cite{HSV2}). Since
$[IT\Rees]=[I\Rees]$, the canonical class of $\widehat\Rees$
vanishes if and only if $[\omega_\Rees]$ is a multiple of
$[I\Rees]$. In view of the theorem this is clearly equivalent to
(i) and (ii) (where in (ii) we now have $c_i=c_j$ for all $i$ and
$j$).
\end{proof}

One should note that the extended Rees algebra (and the associated
graded ring) can have non-trivial projective rank $1$ modules,
even if $\Pic(R)=0$. (Therefore it is not possible to replace
condition (ii) by the requirement that $c_i=c_j$ for \emph{all}
$i$ and $j$.) For example, let $I$ be the intersection of two
maximal ideals in $R=K[X_1,X_2]$. Then it is easy to check by
localization at the maximal ideals of $R$ that the extended Rees
algebra is locally factorial. On the other hand it has divisor
class group isomorphic to $\ZZ$, and all its divisorial ideals are
projective modules.

Another interesting algebra that can be accessed through the Rees
algebra is the \emph{fiber cone} $\Rees/\mm\Rees$, especially in
the situation in which $R=K[X_1,\dots,X_n]$ and
$\mm=(X_1,\dots,X_n)$. If $I$ has a system of generators
$f_1,\dots,f_m$ of constant degree, then one has a natural
embedding $K[f_1,\dots,f_m]$ into $\Rees$. Moreover,
$K[f_1,\dots,f_m]$ is isomorphic to $\Rees/\mm\Rees$, and
therefore a retract of $\Rees$ (see \cite[(2.2)]{Br}). We refer
the reader to the next section where some interesting examples
will be discussed.

\section{Applications}
As mentioned in the introduction, Theorem \ref{main} was inspired
by the its special case derived in \cite{BC2}. We will now use it
in order to extend the results of \cite{BC2} to other classes of
determinantal ideals, and in particular to algebras generated by
minors.

Let $X$ be one of the following types of matrices over a field
$K$:
\begin{itemize}
\item[(G)] an $m\times n$ matrix of indeterminates;

\item[(S)] an $n\times n$ symmetric matrix of indeterminates;

\item[(A)] an $n\times n$ alternating matrix of indeterminates.
\end{itemize}
By $M_t$ we denote the set of $t$-minors of $X$ in the cases (G)
and (S) and the set of $2t$-pfaffians in case (A) (the
$2t$-pfaffians are also elements of degree $t$). In view of the
very detailed analysis of case (G) in \cite{BC2} we restrict
ourselves to an outline containing all the main steps.

There seems to be no single source providing simultaneously all
the details of the cases (G), (S), and (A). For (G) they can be
found in \cite{BC2} and \cite{BV}, for (S) in Abeasis \cite{Ab},
and for (G) in Abeasis and Del Fra \cite{AD} and De Negri
\cite{DN}. We will freely use these sources.

Let $R=K[X]$, $I_t$ be the ideal in $R$ generated by $M_t$, and
$A_t$ the $K$-subalgebra generated by $M_t$. Set
$\Rees_t=\Rees(I_t)$. Since the elements of $I_t$ have all the
same degree, one has a retract
$$
A_t\to \Rees_t\to A_t,
$$
where the embedding $A_t\to \Rees_t$ is induced by the assignment
$f\mapsto fT$, $f\in M_t$, and the kernel of $\Rees_t\to A_t$ is
$\mm R$, with $\mm$ denoting the irrelevant maximal ideal of $R$.
In fact, the bigrading on $\Rees_t$ induces a splitting
$$
\Rees_t=A_t\oplus\mm\Rees_t.
$$
(see \cite[(2.2)]{Br}). It follows that $\mm\Rees_t$ is a prime
ideal.

We assume that the characteristic of $K$ is
\emph{non-exceptional}, i.~e.\ $\chara K=0$ or $\chara
K>\min(t,m-t,n-t)$ in case (G), $\chara K>\min(t,n-t)$ in case
(S), and $\chara K>\min(2t,\allowbreak n-2t)$ in case (A). Then
\begin{equation}
I_t^k=\bigcap_{i=1}^t I_i^{(t-i+1)}.\label{Itk}
\end{equation}

Moreover, if $t<\min(m,n)$, $t<n$, or $2t<n-1$, respectively, the
intersection is irredundant for $k\gg0$; it follows that the
irrelevant maximal ideal is the center of a Rees valuation, and in
particular $\dim A_t=\dim\Rees_t/\mm\Rees_t=\dim R$. In the other
cases the symbolic and the ordinary powers of $I_t$ coincide (and
the canonical module of $\Rees_t$ has been discussed in Bruns,
Simis and Trung \cite{BST}).

That $\Rees_t$ and $A_t$ are Cohen--Macaulay in characteristic $0$
has been shown for all three types in \cite{Br}. In positive
non-exceptional characteristic one finds this result for (G) in
\cite{BC1}, for (A) in Bae\c{t}ica \cite{Ba}, and for (S) in
Bruns, R\"omer and Wiebe \cite{BRW}.

As soon as $\htt I_t>1$, and this is equivalent to $I_t$ being
non-principal, Theorem \ref{main} yields the canonical class of
$\Rees_t$. The hypothesis on the Rees valuations is satisfied in
view of Proposition \ref{PriDec}.

In certain cases the structure of $A_t$ is very easily determined
or classically known:
\begin{enumerate}
\item If $t=1$, then $A_t=R$.

\item If $t=m-1=n-1$ in case (G), $t=n-1$ in case (S), or $2t=n-2$
in case (A), then $A_t$ is isomorphic to a polynomial ring over
$K$. This is easily shown by comparing the Krull dimension with
the number of generators.

\item If $t=n$ in case (S) or $2t=n$ in case (G), then $A_t$ is
isomorphic to a polynomial ring over $K$ for trivial reasons.

\item If $t=\min(m,n)$ in case (G), then $A_t$ is the homogeneous
coordinate ring of a Grassmannian, a factorial Cohen--Macaulay
domain.

\item If $2t=n-1$ in case (A), then $A_t$ is isomorphic to a
polynomial ring over $K$ (this was observed by De Negri and
follows from a theorem of Huneke \cite{Hu}).
\end{enumerate}
\emph{In the following we exclude all these cases, in which the
canonical class is well-known (and trivial).}

The Veronese subring $R^{(t)}$ can be embedded into $R[T]$ in the
same way as $A_t$ into $\Rees_t$, namely by the assignment
$f\mapsto fT$ for all monomials of degree $t$. Then, inside
$\Rees_t$, one obviously has
$$
A_t=\Rees_t\cap R^{(t)}.
$$
Lemma \ref{ReesVal} immediately yields a representation of $A_t$
as an intersection of $R^{(t)}$ with discrete valuations rings. As
in the case (G) treated in \cite{BC2}, one always has the somewhat
surprising equation
$$
A_t=R^{(t)}\cap V_2;
$$
furthermore $R^{(t)}$ is the subintersection of $A_t$ obtained by
omitting the discrete valuation ring $V_2$. This yields an exact
sequence
$$
0\to\ZZ[\pp]\to\Cl(A_t)\to\Cl(R^{(t)})\to0
$$
in which $\pp=P_2\cap A_t$. The group $\Cl(R^{(t)})$ is cyclic of
order $t$. It is not hard to show that $\Cl(A_t)$ is cyclic (of
rank $1$), generated by the class $[\qq]$ of the prime ideal
$fS[T]\cap A_t$, $f$ an arbitrary element of $M_{t+1}$. Moreover,
$[\pp]=-t[\qq]$.

We define an element $\DD\in R$ as follows:
\begin{itemize}
\item[(G)] $\DD$ is the product of all minors of $X$ whose main
diagonals are the parallels to $X_{11},\dots,X_{mm}$. (Such a
parallel starts in each of the $X_{i1}$ and $X_{1j}$.)

\item[(S)] $\DD$ is the product of all minors of $X$, whose main
diagonal is the parallel to $X_{11},\dots,X_{nn}$ starting in one
of the entries $X_{i1}$.

\item[(A)] $\DD$ is the product of all pfaffians whose
anti-diagonal (in the lower triangular part of $X$) is a parallel
to the anti-diagonal of $X$.
\end{itemize}
The valuation associated with $I_j$ is always the (extension of)
the function $\gamma_j$, as defined in De Concini, Eisenbud and
Procesi \cite{DEP}, \cite{Ab}, or \cite{AD}, respectively. For an
element $\delta$ of $M_i$ one has
$$
\gamma_j(\delta)=\begin{cases} 0,& i<j,\\i-j+1,& i\ge j.
\end{cases}
$$
It is now an easy combinatorial exercise to compute
$\gamma_j(\DD)$ in all cases, and it turns out that
$$
\gamma_j(\DD)=\htt I_j=
\begin{cases} (m-j+1)(n-j+1),&\text{(G),}\\
\binom{n-j+1}{2},&\text{(S),}\\
\binom{n-j}{2},&\text{(A)}.
\end{cases}
$$
Thus we have found elements satisfying the condition discussed in
Corollary \ref{Graded}, and
$$
\omega_{\Rees_t}=\DD TR[T]\cap P_1\cap\dots\cap P_t.
$$
We have the presentation $A_t=\Rees_t/\mm \Rees_t$. Now, since
$\mm\Rees_t$ is contained in $P_1$, it follows that
$P_1=\mm\Rees_t$. The (graded) canonical module of $A_t$ is given
by the exact sequence
$$
0\to \Hom_{\Rees_t}(\Rees_t,\omega_{\Rees_t}) \to
\Hom_{\Rees_t}(P_1,\omega_{\Rees_t})\to
\Ext_1^{\Rees_t}(A_t,\omega_{\Rees_t})\to 0.
$$
Let $J$ be the middle term. By divisorial calculation one has
$$
J=\omega_{\Rees_t}:P_1= \DD TR[T]\cap P_2\cap\dots\cap P_t.
$$
and the image of
$\Hom_{\Rees_t}(\Rees_t,\omega_{\Rees_t})=\omega_{\Rees_t}$ is
just $J\cap P_1$. So $\omega_{A_t}=J/J\cap P_1$. In other words,
$\omega_{A_t}$ is the image of $J$ under the epimorphism
$\Rees_t\to A_t$ with kernel $P_1$. Since all our ideals are
bigraded with respect bo the ordinary total degree in R and degree
in $T$, we can replace the image with the intersection:
\begin{align*}
\omega_{A_t}&= \DD TR[T]\cap P_2\cap\dots\cap P_t\cap A_t\\
&=\DD TR[T]\cap P_2\cap A_t.
\end{align*}
The second equation follows since $P_3,\dots,P_t$ meet $A_t$ in
prime ideals of height $>1$. They are superfluous in the
representation of a divisorial ideal.

As remarked above, $P_2\cap A_t=\pp$ has class $-t[\qq]$, and one
can also compute the other term, splitting $\DD T$ into its
factors. The intersection $TR[T]\cap A_t$ is the irrelevant
maximal ideal of $A_t$ and can be omitted, but $fR[T]\cap A_t$,
$f\in M_j$, has class $(j-t)[\qq]$. For the case (G) this has been
computed in \cite[5.3]{BC2}, and the other cases are completely
analogous. A careful count (for (A) one should distinguish the
cases $n$ odd and $n$ even) yields:

\begin{theorem}\label{ClAt}
$\omega_{A_t}=w[\qq]$ with
$$
w=\begin{cases} mn-t(m+n),&\text{\emph{(G),}}\\
\binom{n+1}2-t(n+1),&\text{\emph{(S),}}\\
\binom{n-1}2-2t(n-1),&\text{\emph{(A).}}
\end{cases}
$$
\end{theorem}

\begin{corollary}\label{GorAt}
Set $m=n$ in the cases \emph{(S)} and \emph{(A)}, $u=t$ in the
cases \emph{(G)} and \emph{(S)}, and $u=2t$ in case \emph{(A)}.
Then $A_t$ is Gorenstein if and only if
$$
\frac 1u=\frac 1m+\frac 1n.
$$
\end{corollary}

The reader should note that we have excluded the cases (1)--(5)
above, in the theorem as well as in the corollary. In case (G) we
have only reproduced the results of \cite{BC2}, and our derivation
of Theorem \ref{ClAt} from Theorem \ref{main} has been indicated
in \cite[5.6]{BC2}. The case of Hankel matrices contained in
\cite{BC2} can be treated in the same manner.

\end{document}